\documentstyle[11pt]{article}
\newcommand{\be}{\begin{equation}}
\newcommand{\ee}{\end{equation}}
\newcommand{\bea}{\begin{eqnarray}}
\newcommand{\eea}{\end{eqnarray}}
\newcommand{\ra}{\rightarrow}

\newcommand{\Lra}{\Longleftrightarrow}

\def\Journal#1#2#3#4{{#1} {\bf #2}, #3 (#4)}

\def\PLB{{\em Phys. Lett.} B}

\def\CMP{\em Commun. Math. Phys.}

\def\LMP{\em Lett. Math. Phys.}

\def\JPA{{\em J. Phys.} A}

\begin{document}
\large
\baselineskip=18pt

\thispagestyle{empty}

\hfill {\small DSF-42/01}

\begin{center}
{\Large \bf   Deforming  maps between\\ [3mm]$ sl(n)$, $ sp(2n)$ \\ [3mm] and  
\\[3mm] $U_{q}(sl(n))$, $U_{q}(sp(2n))$} \\[6mm]
 A. Sciarrino\\[3mm]
Universit\`a di Napoli ``Federico II'' \\
Dipartimento di Scienze Fisiche \\
and \\
I.N.F.N. - Sezione di Napoli \\
I-80126 Napoli - Italy
\end{center}

\bigskip

\bigskip

\bigskip

\bigskip

 \centerline{\bf Abstract}

\begin{quote}
Using crystal basis, in the space of symmetric irreducible 
representations, we explicitly   write invertible deforming functionals
between the $q$-deformed universal enveloping algebra and the Lie algebras
$sl(n)$ and $sp(2n)$. For $sl(2)$ we obtain the Curtright-Zachos
map. 
\end{quote}

 \vfill
 {\small \bf
\begin{tabular}{l}
~$\,$Postal adress:  Complesso Universitario di Monte S. Angelo \\ ~$\,$Via 
Cintia - I-80126 Napoli (Italy)\\
 ~$\,$E-mail: SCIARRINO@NA.INFN.IT  
 \end{tabular} }

\clearpage
\thispagestyle{empty}
\mbox{}
\clearpage
\pagestyle{plain}
\setcounter{page}{1}

\section{Introduction}

 Quantum algebras $G_q$ or $U_q (G)$, i.e. the $q$-deformed universal 
enveloping algebra of a semi-simple Lie algebra $G$ \cite{J}, \cite{D}
 admit irreducible highest weight representations (IR) 
 which are for $q$ generic, i.e. not a root of unity, in one-to-one
 correspondence with the IR of the Lie algebra $G$ \cite{L}, \cite{R}. 
 This property strongly
 suggests the existence of a correspondence between the generators of $G$
 and $U_q (G)$. Indeed Curtright and Zachos \cite{CZ} have found an invertible 
 deforming functional ${\cal Q} $  which allows the connection 
 $sl(2) \Lra sl_{q}(2)$.
Denoting by small (resp. capital) letter $j_{\pm,0}$ ($J_{\pm,0}$) the generator
 of $sl(2)$  ($sl_{q}(2)$) it is possible to write ($q$ real number)
 \be
J_{+} = {\cal Q}(j_{\pm},j_{0}) \, j_{+}  \;\;\;\;\;\; 
J_{-} = j_{-} \, {\cal Q}(j_{\pm},j_{0})  \;\;\;\;\;\; 
J_{0} = {\cal Q}_{0}(j_{0}) \, j_{0} =  j_{0}
\ee 
where
\be
{\cal Q} = \sqrt{ \frac{ Ê\{[J_{0} + \mbox{\bf J}]_{q} 
[J_{0} - \mbox{\bf J} - 1]_{q}  \} }
{ Ê\{(j_{0} +  \mbox{\bf j}) (j_{0} -  \mbox{\bf j} - 1)  \} } }
\label{eq:CZ}
\ee
and the operator {\bf j}  ({\bf J})   is defined by the Casimir operator of 
$sl(2)$  ($sl_{q}(2)$) 
\be
C =  \mbox{\bf j} \, ( \mbox{\bf j}+ 1) \;\;\;\;\;\;\;\;\;\;\;\;  
\left (C_{q} = [ \mbox{\bf J}]_{q} \, [ \mbox{\bf J} + 1]_{q} \right)
\ee
In \cite{CZ} it has been argued that the  construction of the invertible 
functional can be generalized for any $sl(n)$, but, at  our knowledge, no proof
has yet been given. It is the aim of this paper to show that indeed such
a functional can explicitly constructed for $sl(n+1)$ and $sp(2n)$  ($n \geq 
1$) in the space of the symmetric IRs, i.e.  the IRs labelled by the 
Dynkin labels $a_{1} > 0,\; a_{i} = 0  \;\; (i > 1)$. In order to prove our statement 
we use the  crystal basis, whose existence
 has been proven for the classical Lie algebra by Kashiwara \cite{Ka}. 
In Sec. 2, in order to make the paper self-contained and to fix the notation,
 we recall the basic definitions of deformed Lie algebra 
 $G_q$ in the Cartan-Chevalley basis and of crystal basis.  In Sec. 3 and Sec. 4 we give our main results. A few
 remarks are given in Sec. 5.

\section{Reminder of deformed algebras and crystal basis}

 Let us recall the definition of $G_q$ associated with a simple Lie algebra
$G$ of rank $r$ defined by the Cartan matrix $ (a_{ij}) $ in the Chevalley
basis. $G_q$ is 
generated by $3r$ elements $e_i^+$, $f_i$ = $e_i^-$ and $h_i$ which satisfy 
$(i,\, j = 1, \ldots , r)$

\be \begin{array}{c}
~[e_i^+,\, e_j^-]  =  \delta_{ij} [h_i]_{q_i} ~~~~~~~~~~~~~~~~ [h_i, h_j] = 0 \\
~~ \\
~[h_i, \, e_j^+]  =  a_{ji} e_j^+ ~~~~~~~~~~~~~~~~~~ [h_i, e_j^-] 
= -a_{ji} e_j^- \label{eq:def}
\end{array}
\ee
where 

\be 
~[x]_q = \frac{q^x - q^{-x}} {q - q^{-1}}  
\ee
and $q_i = q^{d_i}$, $d_i$ being non-zero integers with greatest 
common divisor equal to one such that $d_i \, a_{ij} = d_j\, a_{ji}$.
For simply laced Lie algebra $d_i = 1$, for $sp_q(2n)$  $d_n = 2$.
Further the generators have to satisfy the Serre relations:
\be
\sum_{0 \leq n \leq 1 - a_{ij}} ~ (-1)^n \left[ \begin{array}{c}
1 - a_{ij} \\ n \end{array} \right]_{q_i} \, (e_i^{\pm})^{1 - a_{ij} - n} \, 
e_j^{\pm} (e_i^{\pm})^n =0   \label{eq:dS}
\ee
where
\be
\left[ \begin{array}{c} 
m \\ n \end{array} \right]_q = \frac{[m]_q!}{[m - n]_q! [n]_q}
\ee
\begin{center}
$[n]_q! = [1]_q \: [2]_q \ldots [n]_q$
\end{center}
  In the following
 the deformation parameter $q$  will be assumed different from the roots of 
 the unity. The algebra $G_q$ is also endowed with a Hopf algebra structure,
 but we shall not discuss here this aspect, although very relevant.  

Let us also recall the definition of $gl(n)_{q}$,  ($k,j = 1,2, \dots, n-1; 
i = 1,2, \dots, n-1,n$):
\be \begin{array}{c}
~[e_k^+,\, e_j^-]  =  \delta_{kj} [N_{k} - N_{k+1}]_{q} ~~~~~~~~~~~~~~~~ 
[N_i, N_k] = 0 \\
~~ \\
~[N_i, \, e_j^{\pm}]  = \pm ( \delta_{i,j} \, -  \, \delta_{i-1,j})  \, 
e_j^{\pm} 
 \end{array} \label{eq:gln}
\ee
The Serre relations are computed using $ a_{i,j} = -(\delta_{i-1,j} + 
\delta_{i,j-1} \,)$. So  $gl(n)_{q}$ can be  considered as  $sl(n)_{q} \oplus 
N_n$ with $h_{i} =  N_{i} - N_{i+1}$.

It has been shown in \cite{Ka} that for any $q$-deformed universal 
enveloping algebra of a classical Lie algebra $G$  of rang $r$, in the limit 
$q \ra 0$,
a canonical basis exists, called {\bf crystal basis}, such that, in any
highest weight $\vec{\Lambda}$ IR of $G$: ($i,j = 1, 2 \ldots r$)
\be
\widehat{h}_{i} \, \psi(\vec{\Lambda}; \vec{\lambda} = \{\lambda_{i} \}) =
  \lambda_{i} \, \psi(\vec{\Lambda};  \vec{\lambda} = \{\lambda_{i} \})  \label{eq:f}
  \ee 
\be
 \widehat{e}^{+}_{i} \, \psi(\vec{\Lambda}; \vec{\Lambda}_{i}) = 
 \widehat{e}^{-}_{i} \, \psi(\vec{\Lambda}; -\vec{\Lambda}_{i}) = 0 \label{eq:a}
\ee
\be
 \widehat{e}^{\pm}_{i} \, \psi(\vec{\Lambda}; \vec{\lambda} = \{\lambda_{i} \}) =
  \psi(\vec{\Lambda};  \vec{\lambda} \pm  \vec{a}_{i})  \label{eq:b}
  \ee
  \be
  \vec{\lambda} \pm \vec{a}_{i} = \{\lambda_{j} \pm a_{ij} \} \label{eq:c}
  \ee 
  where we have denoted by $ \widehat{e}^{\pm}_{i}$, $ \widehat{h}_{i} =  h_{i}$
   the generators of 
    $G_{q \ra 0}$, by $\psi(\vec{\Lambda}; \vec{\lambda})$  a generic weight
  $\vec{\lambda}$ state, in the space of the IR with highest weight 
  $ \vec{\Lambda} $, and 
  by $\psi(\vec{\Lambda}; \pm \vec{\Lambda}_{i})$  a state annihilated by
$ \widehat{e}^{\pm}_{i}$. In
\cite{Ka} the relation between $ \widehat{e}^{\pm}_{i}$ and  $e^{\pm}_{i}$ is
 given.
The previous equations imply
  \be
 \widehat{e}^{+}_{i} \,  \widehat{e}^{-}_{j} \, \psi(\vec{\Lambda}; \vec{\lambda}) = 
 \widehat{e}^{-}_{j} \,   \widehat{e}^{+}_{i} \, \psi(\vec{\Lambda}; \vec{\lambda}) = 
\psi(\vec{\Lambda};  \vec{\lambda} + \vec{a}_{i} - \vec{a}_{j} ) \label{eq:d}
 \ee
which for $ i = j$ can be written
  \bea
&  \widehat{e}^{+}_{i} \,  \widehat{e}^{-}_{i} \, \psi(\vec{\Lambda}; \vec{\lambda}) = 
\mbox{\bf 1} \, \psi(\vec{\Lambda};  \vec{\lambda}) \;\;\;\;\;\; (\vec{\lambda} 
\neq  - \vec{\Lambda_{i}})
\nonumber \\
&  \widehat{e}^{-}_{i} \,  \widehat{e}^{+}_{i} \, \psi(\vec{\Lambda}; \vec{\lambda}) = 
 \mbox{\bf 1} \, \psi(\vec{\Lambda};  \vec{\lambda}) \;\;\;\;\;\; (\vec{\lambda} 
\neq   \vec{\Lambda_{i}})
 \label{eq:e}
 \eea
 Note that from eq.(\ref {eq:d}) the Serre relations for $G$ are trivially
 satisfied by $\{ \widehat{e}^{\pm}_{i} \}$.

\section{Relation between $sl_{q}(n)$ and $sl(n)$}
 
 In this section we derive explicitly the invertible functionals which
 connect $sl_{q}(n)$ and $sl(n)$. 
 Let us recall, e.g. see \cite{FSS}, that a state of an IR of $sl(n)$ is 
 identified by a Young tableaux, which  is a  pattern of $n$ rows, the i-th row
 containing $l_{i}$ boxes, where $l_{i}$ are not negative integer 
 satisfying, for the highest weight state, $l_{k} \geq l_{k+1}$ 
 ($k = 1,2,\ldots,n-1$).  Strictly speaking one should consider class of
 Young tableaux as two tableaux differing for a block of ($ p \times n$) 
 boxes i.e. by  the first $p$ columns ($ 1 \leq p \leq l_{n}$), identify the same 
 IR. The relation 
 between the $n-1$ Dynkin label $a_{k}$ and the $n$ labels $l_{i}$ is
 \be
 a_{k} = l_{k} - l_{k+1}
 \ee 
 The integer $l_{i}$ is the eigenvalue of the operator $N_{i}$ on the
 state identified by the corresponding Young tableaux and 
 $\lambda_{i} =  l_{i} -  l_{i+1}$. 
   From eq.(\ref{eq:gln}) it follows
  \be
 ~[N_i, \,  \widehat{e}_k^{\pm}]  = \pm ( \delta_{i,k} \, -  \, \delta_{i-1,k})  \, 
 \widehat{e}_k^{\pm}  \label{eq:numero}§§
\ee   
 As first step we prove the following proposition
 
{\bf Prop. 1} - The generators ($i = 1,2,\ldots,n-1$)
 \be
 E^{+}_{i} =  \widehat{e}^{+}_{i} \, \sqrt{(N_{i} \, + \, 1) \, N_{i+1}} 
 \;\;\;\;\;\;
  E^{-}_{i} = \sqrt{(N_{i} \, + \, 1) \, N_{i+1}} \,\,\,  \widehat{e}^{-}_{i} 
  \;\;\;\;\;\; 
 H_{i} =  N_{i} \, - \, N_{i+1} \label{eq:sl}
\ee
where $\{ \widehat{e}^{\pm}_{i},  N_{i}, N_{n} \}$ are the   generators in 
the crystal basis of $gl(n)_{q \ra 0}$, satisfy the defining relations, in 
the Cartan-Chevalley basis, of $sl(n)$ in the space of the symmetric IRs
labelled by the Young tableaux $l_{1} = \Lambda,\; l_{i} = 0$.

Proof: On the highest weight state we have
\bea
& [ E^{+}_{i}, \,  E^{-}_{i}] \,\, \psi(\vec{\Lambda};  \vec{\Lambda}) =
  E^{+}_{i} \,  E^{-}_{i} \,\, \psi(\vec{\Lambda}; \vec{\Lambda}) = 
   \nonumber \\
&  N_{i} \,  (N_{i+1} \, +  \, 1)  \,\, \psi(\vec{\Lambda}; \vec{\Lambda}) = 
  (l_{i}l_{i+1} \,  + \,  l_{i} §)  \,\, \psi(\vec{\Lambda}; \vec{\Lambda})
\label{eq:hw}
\eea 
The eigenvalue of the r.h.s. of eq.(\ref{eq:hw}) has to identified with
the eigenvalue of $H_{i}$, i.e. $l_{i} - l_{i+1}$, which implies 
$l_{1} = \Lambda, \; l_{i} = 0 \;\; (i \geq 1)$.
Then we  consider the states $\{\psi(\vec{\Lambda}; - 
\vec{\Lambda}_{i})\}$ which are eigenstates with vanishing eigenvalue of $N_{i}$
\be   
N_{i} \,\, \psi(\vec{\Lambda}; - \vec{\Lambda}_{i}) = 0 \label{eq:van}
\ee
 We have, from eq.(\ref{eq:sl}):
\bea
& [ E^{+}_{i}, \,  E^{-}_{i}]  \,\, \psi(\vec{\Lambda};  
- \vec{\Lambda}_{i}) = \nonumber \\
& - E^{-}_{i} \,  E^{+}_{i} \,\, \psi(\vec{\Lambda}; - \vec{\Lambda}_{i}) = 
- (N_{i} \, + \, 1) \,N_{i+1} \,\, \psi(\vec{\Lambda}; - \vec{\Lambda}_{i}) = 
 \nonumber \\
 & - N_{i+1} \,\, \psi(\vec{\Lambda}; - \vec{\Lambda}_{i}) =
  H_{i}  \, \, \psi(\vec{\Lambda}; - \vec{\Lambda}_{i})
\eea  
On the state $\psi(\vec{\Lambda}; \vec{\lambda})$ 
($\vec{\lambda} \neq  - \vec{\Lambda}_{i}$), from 
eqs.(\ref{eq:b})-(\ref{eq:e}) we have
\be
[H_{i}, \,   E^{\pm}_{j}]\,\,\psi(\vec{\Lambda}; \vec{\lambda}) = 
\pm \, a_{ji} \,  E^{\pm}_{j}  \,\, \psi(\vec{\Lambda}; \vec{\lambda})  
\ee 
  \be
  [ E^{+}_{i}, \,  E^{-}_{j}] \,\, \psi(\vec{\Lambda}; \vec{\lambda}) = 
  \delta_{ij} \, H_{i} \,\, \psi(\vec{\Lambda}; \vec{\lambda})  
  \ee 
So we can write
  \be
[H_{i}, \,   E^{\pm}_{j}]_{\psi} = \pm \, a_{ji} \,  E^{\pm}_{j} \label{eq:L1}
\ee 
  \be
  [ E^{+}_{i}, \,  E^{-}_{j}]_{\psi}  = \delta_{ij} \, H_{i} \label{eq:L2}
  \ee
  where the lower label $\psi$  in the commutator reminds that the 
  relations hold when applied on a state of a symmetric IR and not
  as general algebraic expressions.
   Finally we have to prove that generators $E^{\pm}_{i}$ satisfy the Serre 
 relation for $sl(n)$:
 \be
\sum_{0 \leq n \leq 1 - a_{ij}} ~ (-1)^n \left( \begin{array}{c}
1 - a_{ij} \\ n \end{array} \right) \, (E_i^{\pm})^{1 - a_{ij} - n} \, 
E_j^{\pm} (E_i^{\pm})^n =0 
\ee
 This is an immediate consequence of eq.(\ref{eq:sl}) and of the following 
 identity 
 ($a \geq 1$, $ a \in \mbox{\boldmath $Z_{+}$}$, 
  $ z \in \mbox{{\bf C}} $), evaluated (for $E_i^{+}$)
 for  $N = N_{i+1} + 1$ , $ q = 1$,   and $a = z = 1$
 \be
\sum_{0 \leq n \leq (1 + a)} ~ (-1)^n \left [ \begin{array}{c}
1 + a \\ n \end{array} \right ]_{q} \, [N - nz]_{q} = 0  \label{eq:id}
\ee
 
Then we easily  get 

{\bf Prop. 2} - The generators ($i = 1,2,\ldots,n-1$)
 \be
 e^{+}_{i} =  \widehat{e}^{+}_{i} \, \sqrt{[N_{i} \, + \, 1]_{q} \,[N_{i+1}]_{q}} 
 \;\;\;\;\;\;
  e^{-}_{i} = \sqrt{[N_{i} \, + \, 1]_{q} \, [N_{i+1}]_{q}} \,\,\,
   \widehat{e}^{-}_{i} 
  \;\;\;\;\;\; 
 h_{i} = H_{i}  \label{eq:slq}
\ee
define $sl_{q}(n)$ in the Cartan-Chevalley basis in the space of the symmetric 
IRs.

Proof: The second defining relation in eq.(\ref{eq:def}) is immediately 
proven while the first one follows by the identity
\be
[N_{i}]_{q} \, [N_{i+1} \, + \, 1]_{q} - [N_{i} \, + \, 1]_{q} \, [N_{i+1}]_{q} =
[N_{i} \, - \, N_{i+1}]_{q}
\ee
The deformed Serre relations  eq.(\ref{eq:dS}) are proven using  
eq.(\ref{eq:id}).
  
From {\bf Prop.1} and {\bf Prop.2}, in the symmetric basis,
 it follows the  relation between $sl(n)$ and $sl_{q}(n)$ 
\be
e^{+}_{i} = E^{+}_{i} \, \sqrt{\frac{[N_{i} \, + \, 1]_{q} \, [N_{i+1}]_{q}}
{(N_{i} \, + \, 1) \, N_{i+1}}}  \;\;\;\;\;\;   
e^{-}_{i} =  \sqrt{\frac{[N_{i} \, + \, 1]_{q} \, [N_{i+1}]_{q}}
{(N_{i} \, + \, 1) \, N_{i+1}}} \, E^{-}_{i}  \label{eq:I}
\ee
Let us discuss in more detail the case $sl(2)$. In this case all the IRs are  of the type 
  we have called  symmetric  and  the eq.(\ref{eq:I}) reads
\be
j^{+} = J^{+}  \, \sqrt{\frac{[N_{1} \, + \, 1]_{q} \, [N_{2}]_{q}}
{(N_{1} \, + \, 1) \, N_{2}}}  \;\;\;\;\;\;   
j^{-}  =  \sqrt{\frac{[N_{1} \, + \, 1]_{q} \, [N_{2}]_{q}}
{(N_{1} \, + \, 1) \, N_{2}}} \, J^{-}  \label{eq:Ia}
\ee 
From the following relation between the operators {\bf j} and $j_{0}$ and 
our operators 
 \be
\mbox{\bf j} = \frac{N_{1} \, + \,  N_{2}}{2}   \;\;\;\;\;\; 
j_{0} =  \frac{N_{1} \, - \,  N_{2}}{2}
\ee
we see that eq.(\ref{eq:Ia}) is equivalent to the Curtright-Zachos formula
eq.(\ref{eq:CZ}). Moreover, in this case, the
 deforming map can be turned into a general algebraic relation. 
 As a first consequence of eq.(\ref{eq:I}),  we can give a new construction of 
 $so_{q}(3)$ in terms of standard bosons. Let us recall that
 Van der Jeugt \cite{Van} has shown that the following generators satisfy
 the commutation of $so_{q}(3)$
 \bea
L_{+} &  =  & q^{N_{-1}} \, q^{-N_{0}/2} \, \sqrt{q^{N_{1}} + q^{-N_{1}}} \, 
b_1^+ b_0 \nonumber \\ 
& +  &  b_0^+ b_{-1} \,  q^{N_{1}} \, q^{-N_{0}/2} \, 
\sqrt{q^{N_{-1}} + q^{-N_{-1}}}   \nonumber \\ 
  L_{-} &  =  &  b_0^+ b_{1} \, q^{N_{-1}} \, q^{-N_{0}/2} \, 
  \sqrt{q^{N_{1}} + q^{-N_{1}}} \nonumber \\Ê 
  & +  & \,  q^{N_{1}} \, q^{-N_{0}/2} \, 
\sqrt{q^{N_{-1}} + q^{-N_{-1}}} \,  b_{-1}^+ b_{0}   \nonumber \\Ê 
    L_{0} & = & N_1  - N_{-1}   
\eea
where $b_{\pm,0},\; b^{+}_{\pm,0}$ are the Biedenharn-MacFarlane \cite{B}, 
\cite{M} $q$-bosons.
From our result it follows that the following generators satisfy, on the
states of any IR, the commutation relations of $so_{q}(3)$
 \bea
L_{+} &  =  & \sqrt{2} \, (\tilde{b}_1^+ \tilde{b}_0  +  \tilde{b}_0^+ 
\tilde{b}_{-1})  \,  \sqrt{\frac{[N_{1} \, + \, 1]_{q} \, [N_{2}]_{q}}
{(N_{1} \, + \, 1) \, N_{2}}}  \nonumber \\
L_{-} &  =  &    \sqrt{\frac{[N_{1} \, + \, 1] \, [N_{2}]_{q}}
{(N_{1} \, + \, 1) \, N_{2}}}  \,\,\, \sqrt{2} \, (\tilde{b}_1 \tilde{b}_0^+  + 
 \tilde{b}_0  \tilde{b}_{-1}^+)  \nonumber \\
 N_{1}  & =  & 2 \; \tilde{b}_1^+ \tilde{b}_1 \, + \, \tilde{b}_0^+ \tilde{b}_0 
 \nonumber \\ 
 N_{2}  & =  & 2 \; \tilde{b}_{-}1^+ \tilde{b}_{-1} \, +  \, \tilde{b}_0^+ \tilde{b}_0
\nonumber \\   L_{0} & = & (N_1  - N_{2} )/2  =  \tilde{b}_1^+ \tilde{b}_1 
\, - \, \tilde{b}_{-}1^+ \tilde{b}_{-1}  
\eea
 where $\tilde{b}_{\pm,0}^+,\; \tilde{b}_{\pm,0}$ are standard bosonic 
 operators. 

\section{ Relation between $sp_{q}(2n)$ and $sp(2n)$}

In this section we derive explicitly the invertible functionals which
 connect $sp_{q}(2n)$ and $sp(2n)$.  Let us recall, e.g. see \cite{FSS}, 
 that   an IR of $sp(2n)$  can be 
 identified by a $n$-rows Young tableaux, the i-th row  
 containing $l_{i}$ boxes, $l_{k} \geq l_{k+1}$ ($k = 1,2,\ldots,n-1$).    
  The relation 
 between the $n$ Dynkin label $a_{k}$ and the $n$ labels $l_{i}$ is
 \be
 a_{k} = l_{k} - l_{k+1} \;\;\;\;\;\;\;\; a_{n} = l_{n}  
 \ee 
 The integer $l_{i}$ is the eigenvalue of the operator $N_{i}$ on the
 state identified by the corresponding Young tableaux. 
 Now we  define 
 \be
 ~[N_n, \,  \widehat{e}_j^{\pm}]  = \pm ( 2 \, \delta_{n,j} \, -  \, \delta_{n-1,j})  \, 
 \widehat{e}_j^{\pm}
\ee   
while the action of for $N_j, j \neq n$, is given by eq.(\ref{eq:numero}). 
 In complete analogy with the previous case, using eq.(\ref{eq:van}), we  get
  
  {\bf Prop. 3} - The generators ($k = 1,2,\ldots,n-1$)
 \bea
& E^{+}_{k} =  \widehat{e}^{+}_{k} \, \sqrt{(N_{k}  +  1) \, N_{k+1}} 
 \;\;\;\;\;\;\;\;
  E^{-}_{k} = \sqrt{(N_{k}  +  1) \, N_{k+1}} \,\,\,  \widehat{e}^{-}_{k} 
 \nonumber \\ 
& E^{+}_{n} = \frac{ \widehat{e}^{+}_{n}}{2} \, \sqrt{(N_{n}  +  1) \, (-  N_{n} \, - 
 \, 2)} 
 \;\;\;\;\;\;\;\;
  E^{-}_{n} = \sqrt{(N_{n} \, + \, 1) \, (- N_{n}  -  2)} \,\,\, 
  \frac{ \widehat{e}^{-}_{n}}{2}  
   \nonumber \\
& H_{k} =  N_{k}  -  N_{k+1} \;\;\;\;\;\;\;\;  H_{n} = N_{n}  + \frac{1}{2}   
\label{eq:sp}
\eea
 define, in the Cartan-Chevalley basis,  $sp(2n)$ in the  spaces of 
 symmetric IRs labelled 
 by the Young tableaux with $l_{1} = \Lambda, \; l_{i} = 0 \;\; (i \geq 1)$.
 
 To prove the Serre relations between  $E^{+}_{n-1}$ and $E^{+}_{n}$ one
 has to use  eq.(\ref{eq:id}) evaluated for $N = N_{n}$,  $q = 1$, $z = -2$ 
 and $ a = 2 $, for $a_{n,n-1}$, and  $ a = 3$, for $a_{n-1,n}$. 
  
  {\bf Prop. 4} - The generators ($ k = 1,2,\ldots,n-1$)
 \bea
&  e^{+}_{k} =  \widehat{e}^{+}_{k} \, \sqrt{[(N_{k} \, + \, 1]_{q} \, [N_{k+1}]_{q}} 
 \;\;\;\;\;\; 
  E^{-}_{k} = \sqrt{[(N_{k} \, + \, 1]_{q} \, [N_{k+1}]_{q}} \,\,\,  \widehat{e}^{-}_{k} 
  \nonumber \\
& e^{+}_{n} =  \frac{1}{q + q ^{-1}}  \, \widehat{e}^{+}_{n} \,
 \sqrt{[ N_{n}  +  1]_{q} \, [-  N_{n}  -  2)]_{q}} 
 \;\;\;\; 
  E^{-}_{n} = \frac{1}{q + q ^{-1}}  \,\sqrt{[(N_{n}  +  1]_{q} \, 
  [- N_{n}  -  2)]_{q} } 
  \,\,\,  \widehat{e}^{-}_{n} 
   \nonumber \\
& h_{k} =  H_{k}  \;\;\;\;\;\;\;\;   h_{n} = H_{n}  \label{eq:spq}
\eea
define $sp_{q}(2n)$ in the Cartan-Chevalley basis in the space of the symmetric 
IRs.

The  first defining relation in eq.(\ref{eq:def}) is  
proven using the identity
\be
[N - 1]_{q} \, [- N]_{q} - [N  + 1]_{q} \, [- N - 2]_{q} =
[2 N + 1 ]_{q} \times (q + q ^{-1})
\ee

 From {\bf Prop.3} and {\bf Prop.4} it follows 
  the  relation between $sp(2n)$ and $sp_{q}(2n)$ 
\bea
& e^{+}_{k} = E^{+}_{k} \, \sqrt{\frac{[(N_{k} \, + \, 1]_{q} \, [N_{k+1}]_{q}}
{(N_{k} \, + \, 1) \, N_{k+1}}}  \;\;\;\;\;\;   
e^{-}_{k} =\sqrt{\frac{[(N_{k} \, + \, 1]_{q} \, [N_{k+1}]_{q}}
{(N_{k} \, + \, 1) \, N_{k+1}}} \,\,\, E^{-}_{k}  \nonumber \\
& e^{+}_{n} =  \frac{2}{q + q ^{-1}}  \, E^{+}_{n} \, \sqrt{\frac{[( N_{n}  +  1]_{q}  \,[- N_{n}  - 
 2)]_{q}}{( N_{n}  +  1) \,(- N_{n}  -  2)}}  \;\;\;\;\;\;   
e^{-}_{n} =   \frac{2}{q + q ^{-1}}  \, \sqrt{\frac{[(N_{n}  +  1]_{q} \,[- N_{n} \, - 
\, 2)]_{q}}{(N_{n}  +  1) \,(- N_{n}  -  2)}} \,\,\, E^{-}_{n}
     \label{eq:II}
\eea

 \section {\bf   Conclusions}
 
 We have found explicit invertible maps between $U_{q}(sl(n))$, 
$U_{q}(sp(2n))$  and $ sl(n)$, $ sp(2n)$ which hold on the states of
the IRs labelled by  the   Dynkin label $a_{1}$, that is the
  symmetric  representations. Let us remark that the Cartan sub-algebra is 
  left undeformed and it can also be identified with the set of diagonal
  generators in the crystal basis. As a byproduct result  we have   
obtained:
\begin{itemize}
\item  in the spaces of the  symmetric IRs, explicit expressions of the
generators of $U_{q \ra 0}(sl(n))$ and  $U_{q \to 0}(sp(2n))$  in terms of the
generators of  $ sl(n)$ and $ sp(2n)$. Note however that the correspondence
$ E^{\pm}_{i} \Lra  \widehat{e}^{\pm}_{i}$  has to be handled with some 
precautions as the 
factor  $(N_{i} \, + \, 1) \, N_{i+1}$ appearing in {\bf Prop.1} and in {\bf Prop.3}
  is vanishing  on  some states.  This is an interesting result
as the relations  given in \cite{Ka}  between $ \widehat{e}^{\pm}_{i}$ and 
 $e^{\pm}_{i}$ are rather cumbersone for Lie algebra $G$ of rank
 larger than one.
\item from eq.(\ref{eq:sl})-(\ref{eq:slq}) 
(resp. eqs.(\ref{eq:sp})-(\ref{eq:spq}))
an expression of the action of  the generators of $sl(2)$ (resp. $sl_{q}(2)$) 
embedded in $sl(n)$ and $sp(2n)$ (resp. $sl_{q}(n)$ and $sl_{q}(2n)$ on the 
states of the symmetric IRs.
\end{itemize}


\begin{thebibliography}{99}

\bibitem{J}  M. Jimbo,{\it ``A $q$-difference analogue of U($g$) and the
Yang-Baxter equation"}, \Journal{\LMP} {10} {63} {1985}

\bibitem{D} V.G. Drinfeld,
 {\it ``Quantum Groups"}, in {\em Proc.Int.Congr. of Math.}, MSRI
 Berkely, California (1986) 

\bibitem{L} G. Lusztig,
 {\it `` Modular representations of Quantum Groups"},  Preprint (1988) 
 
 \bibitem{R}  M.  Rosso, {\it `` Finite dimensional representations of the
 quantum analog the algebra of a complex simple Lie algebra"}, 
 \Journal{\CMP }{117}{581}{1988}
   
\bibitem{CZ} T.L. Curtright and  C.K. Zachos, {\it ``Deforming maps for  
quantum algebras"}, \Journal{\PLB} {243} {237} {1990}

\bibitem{Ka}  M. Kashiwara, {\it `` Crystalizing the $q$-Analogue   
of Universal Enveloping Algebras"}, \Journal{\CMP }{133}{249}{1990}

\bibitem{FSS} L. Frappat, A. Sciarrino and P. Sorba, {\it ``Dictionary on 
Lie Algebras and Superalgebra"}, Academic Press, London (2000)

\bibitem{Van} J. Van der Jeugt, {\it `` On the principal subalgebra of quantum
enveloping algebras $gl_{q}(l + 1)$"}, \Journal{\JPA} {25} {L213} {1992} 

\bibitem{B}  L.C. Biedenharn,{\it ``The quantum group $SU_{q}(2)$ and 
a $q$-analogue of the boson operators"}, \Journal{\JPA} {22} {L873} {1989}

\bibitem{M}  A.J. Macfarlane, {\it ``On $q$-analogues of the quantum harmonic 
oscillator and the  quantum group $SU_{q}(2)$"}, \Journal{\JPA} {22} {4581} {1989}

\end{thebibliography}
\end{document}